\pdfoutput=1
\RequirePackage{ifpdf}
\ifpdf 
\documentclass[pdftex]{sigma}
\else
\documentclass{sigma}
\fi

\numberwithin{equation}{section}
\newtheorem{Theorem}{Theorem}[section]
\newtheorem{Lemma}[Theorem]{Lemma}
\newtheorem{Proposition}[Theorem]{Proposition}
{\theoremstyle{definition}
\newtheorem{Definition}[Theorem]{Definition}
\newtheorem{Remark}[Theorem]{Remark}
}

\begin{document}

\newcommand{\arXivNumber}{1402.0072}

\allowdisplaybreaks

\renewcommand{\thefootnote}{$\star$}

\renewcommand{\PaperNumber}{057}

\FirstPageHeading

\ShortArticleName{Induced Representations and Hypergroupoids}

\ArticleName{Induced Representations and Hypergroupoids\footnote{This paper is a~contribution to the Special Issue on
Noncommutative Geometry and Quantum Groups in honor of Marc A.~Rief\/fel.
The full collection is available at
\href{http://www.emis.de/journals/SIGMA/Rieffel.html}{http://www.emis.de/journals/SIGMA/Rieffel.html}}}

\Author{Jean RENAULT}

\AuthorNameForHeading{J.~Renault}

\Address{Universit\'e d'Orl\'eans et CNRS (UMR 7349 et FR2964), D\'epartement de Math\'ematiques,\\
F-45067 Orl\'eans Cedex 2, France} 
\Email{\href{mailto:Jean.Renault@univ-orleans.fr}{Jean.Renault@univ-orleans.fr}}
\URLaddress{\url{http://www.univ-orleans.fr/mapmo/membres/renault/}}

\ArticleDates{Received February 01, 2014, in f\/inal form May 26, 2014; Published online June 03, 2014}

\Abstract{We review various notions of correspondences for locally compact groupoids with Haar systems, in particular
a~recent def\/inition due to R.D.~Holkar.
We give the construction of the representations induced by such a~correspondence.
Finally, we extend the construction of induced representations to hypergroupoids.}

\Keywords{groupoids; $C^*$-algebras; correspondences; induced representations; hypergroups}

\Classification{22D30; 22D25; 22A22; 46L08}

\rightline{\it Dedicated to Professor Marc A.~Rieffel on the occasion of his 75th birthday}

\renewcommand{\thefootnote}{\arabic{footnote}}
\setcounter{footnote}{0}

\section{Introduction}

\looseness=-1
Among many great accomplishments, M.~Rief\/fel is well known for the theory of strong Morita equivalence of $C^*$-algebras
he introduced in his inf\/luential article~\cite{rie:induced} on induced representations for $C^*$-algebras.
He develops there a~$C^*$-algebraic setting which generalizes much of G.~Mackey's theory of induced representations for
locally compact groups.
As described in~\cite{rie:kingston}, transformation group $C^*$-algebras provide many examples of strongly Morita
equivalent $C^*$-algebras.
The notion of groupoid equivalence introduced in~\cite{ren:halifax, ren:kingston} to give a~common framework to these
examples is directly inspired by Rief\/fel's theory.
The key fact (stated in~\cite[Section~3]{ren:kingston} and proved in~\cite[Theorem~2.8]{mrw:Morita}) is that a~groupoid
equivalence implements a~strong Morita equivalence of the groupoid $C^*$-algebras.
However, concerning induced representations, M.~Rief\/fel points out that~\cite{rie:induced} does not cover some parts of
Mackey's theory such as the subgroup theorem or the intertwining number theorem.
This raises the question to f\/ind a~proper setting for these theorems.
For that purpose, it seems that a~necessary step is to f\/ind an analogue at the groupoid level of the notion of
$C^*$-correspondence.
This $C^*$-algebraic notion of correspondence appears naturally in~\cite[Definition~4.19]{rie:induced} under the name of
Hermitian~$B$-rigged~$A$-module (we shall simply say $(A,B)$-$C^*$-correspondence, or $C^*$-correspondence from~$A$
to~$B$).
They (or rather their isomorphism classes) can be viewed as morphisms in a~category having $C^*$-algebras as objects;
Morita equivalences are isomorphisms in this category.
Since~\cite{rie:induced}, various notions of correspondences for groupoids have been proposed.
They will be reviewed in Section~\ref{section2} of this article.
However, it is only recently that, building upon the previous def\/initions, a~satisfactory notion has been found by
R.D.~Holkar.
It has two big advantages: f\/irst, it makes the groupoid $C^*$-algebra construction a~functor from the category of
groupoids to the category of $C^*$-algebras and second, it implements many classical cases of induced representations.
In Section~\ref{section3}, we shall describe the construction of representations induced through a~groupoid correspondence and give
some examples which illustrate this notion.
In Section~\ref{section4}, we shall see the necessity to enlarge the ca\-te\-go\-ry of groupoids to capture some natural constructions of
induced representations.
This will lead us to the notion of locally compact hypergroupoid with Haar system and to the construction of their
$C^*$-algebras.
An important class of examples is provided by double coset hypergroups, and more generally by spatial hypergroupoids
which have been considered earlier in~\cite{hr:2013}.

The groupoid $C^*$-algebra construction works very well for non-Hausdorf\/f group\-oids (more precisely topological
groupoids such that each point has a~compact Hausdorf\/f neighborhood, the unit space is Hausdorf\/f and there is
a~continuous Haar system).
Nevertheless, for the sake of simplicity, we shall assume here that our locally compact spaces and groupoids are
Hausdorf\/f and second countable.
Most of the theory goes through in the non-Hausdorf\/f case but a~systematic treatment would impair the legibility of the
exposition.
When we consider Borel spaces and groupoids, we implicitly assume that they are analytic.

We use the terminology and the notation of~\cite{adr:amenable}.
The unit space of a~groupoid~$G$ is denoted by~$G^{(0)}$.
The elements of~$G$ are usually denoted by $\gamma, \gamma',\ldots$; those of $G^{(0)}$ are denoted by $x,y,\ldots$ or
$u,v\ldots$.
The structure of~$G$ is def\/ined by the inclusion map $i: G^{(0)}\rightarrow G$ (we shall identify~$x$ and $i(x)$),
the range and source maps $r,s:G\rightarrow G^{(0)}$,
the inverse map $\gamma\mapsto\gamma^{-1}$ from~$G$ to~$G$ and the
multiplication map $(\gamma,\gamma')\mapsto \gamma\gamma'$ from the set of composable pairs
\begin{gather*}
G^{(2)}=\{(\gamma,\gamma')\in G\times G: s(\gamma)=r(\gamma')\}
\end{gather*}
to~$G$.
Given $A,B\subset G^{(0)}$, we write $G^A=r^{-1}(A)$, $G_B=s^{-1}(B)$ and $G^A_B=G^A\cap G_B$.
Similarly, given $x,y\in G^{(0)}$, we write $G^x=r^{-1}(x)$, $G_y=s^{-1}(y)$ and $G(x)=G^x_x$.
We assume that $G^{(0)}\subset G$ and $G^{(2)}\subset G\times G$ have the subspace topology.
We also include in the def\/inition of a~topological groupoid the assumptions that the range and source maps are open.

Most spaces occuring in the theory of groupoids are f\/ibered spaces: we shall say that a~set~$Y$ is {\it fibered} over
a~set~$X$ if a~surjective map $p:Y\rightarrow X$, called the {\it projection map}, has been specif\/ied.
When two sets~$Y$ and~$Z$ are f\/ibered over~$X$ via the maps $p:Y\rightarrow X$ and $q:Z\rightarrow X$, we def\/ine the
{\it fibered product}
\begin{gather*}
Y*Z=\{(y,z)\in Y\times Z: p(y)=q(z)\}.
\end{gather*}
It is also a~set f\/ibered over~$X$.
A~{\it left~$G$-space} consists of a~set~$X$ f\/ibered over $G^{(0)}$ by a~map $r_X:X\rightarrow G^{(0)}$, called the
projection or {\it moment map} and a~map $(\gamma,x)\in G*X\mapsto \gamma x\in X$, called the {\it action map},
where~$G$ is f\/ibered over $G^{(0)}$ by the source map, such that the following equalities hold whenever they make sense:
$\gamma(\gamma'x)=(\gamma\gamma')x$ and $ux=x$ for a~unit $u\in G^{(0)}$.
The relation $x\sim y$ if and only if there exists $\gamma\in G$ such that $x=\gamma y$ is an equivalence relation and
the quotient space is denoted by $G\backslash X$.
One def\/ines similarly a~{\it right~$G$-space~$X$}.
Then the moment map is denoted by $s_X$ and the quotient space by $X/G$.
An {\it $(H,G)$-space}, where $H$, $G$ are groupoids, is a~space~$X$ which is both a~left~$H$-space and a~right~$G$-space
and such that $(hx)\gamma=h(x\gamma)$ for all $(h,x,\gamma)\in H*X*G$.
In the locally compact setting, we assume that~$X$ and~$G$ are locally compact, the moment map is continuous and open
and the action map is continuous.
One says that the left~$G$-space~$X$ is {\it free} if the map from $G*X$ to $X\times X$ sending $(\gamma,x)$ to $(\gamma
x,x)$ is one-to-one and that it is {\it proper} if this map is proper.
One says that the groupoid~$G$ is {\it proper} if the left~$G$-space $G^{(0)}$ is proper.
In the Borel setting, we assume that~$X$ and~$G$ are Borel and that the moment map and the action map are Borel.
Given a~left~$G$-space~$X$, one can form the {\it semi-direct product groupoid} $G\ltimes X$ (see, e.g.,~\cite[Subsection~2.1.a]{adr:amenable}): as a~set, it is $G*X$ but it is sometimes preferable to write its elements as
$(\gamma x, \gamma, x)$ rather than $(\gamma, x)\in G*X$ to let appear the range of the element.
Its unit space is~$X$ and $i(x)=(x,r_X(x),x)$; $r (\gamma x, \gamma, x)=\gamma x$ and $s(\gamma x, \gamma, x)=x$; its
multiplication is given by
\begin{gather*}
(\gamma' \gamma x, \gamma', \gamma x)(\gamma x, \gamma, x)=(\gamma'\gamma x, \gamma'\gamma, x)
\end{gather*}
and its inverse map by $(\gamma x, \gamma, x)^{-1}= (x, \gamma^{-1}, \gamma x)$.
In the locally compact setting, $G\ltimes X$ becomes a~locally compact groupoid.
One def\/ines similarly the semi-direct product groupoid $X\rtimes G$ for a~right~$G$-space~$X$.

\looseness=-1
Given locally compact (resp.~Borel) spaces $X$, $Y$ and a~continuous (resp.~Borel) projection map $p:Y\rightarrow X$, a~continuous (resp.~Borel) {\it~$p$-system of measures} is a~family $\alpha=(\alpha_x)_{x\in X}$ such that $\alpha_x$ is a~Radon
[resp.~$\sigma$-f\/inite] measure on $p^{-1}(x)$ and such that for all function~$f$ on~$Y$ complex-valued, continuous with
compact support (resp.\ positive and Borel), $\displaystyle
x\mapsto \int f d\alpha_x$ is continuous (resp.\
Borel) on~$Y$; this function is denoted $\alpha(f)$.
If $X$, $Y$ are (say right)~$G$-spaces and~$p$ is equivariant, we say that~$\alpha$ is equivariant if for all $\gamma\in
G$, $\alpha_{r(\gamma)}\gamma=\alpha_{s(\gamma)}$.
A~{\it $($left$)$ Haar system} for the locally compact groupoid~$G$ is a~continuous and equivariant~$r$-system of measures
$\lambda=(\lambda^x)_{x\in G^{(0)}}$ for the left~$G$-space~$G$.
It is an essential piece of data for def\/ining the convolution algebra $C_c(G,\lambda)$ and the reduced and full
$C^*$-algebras $C^*_r(G,\lambda)$ and $C^*(G,\lambda)$ (which are usually denoted simply by $C^*_r(G)$ and $C^*(G)$).
When~$G$ is an \'etale groupoid, i.e.~the range map~$r$ is a~local homeomorphism, one chooses implicitly the Haar system
$\lambda=(\lambda^x)_{x\in G^{(0)}}$ where $\lambda^x$ is the counting measure on $G^x$.

\section{Groupoid correspondences}\label{section2}

\subsection{Previous def\/initions}\label{section2.1}

The notion of groupoid correspondence was introduced at the same time as the notion of groupoid equivalence.
In~\cite{ren:halifax}, a~{\it correspondence} from a~locally compact groupoid~$H$ to another~$G$ is def\/ined as
an~$(H,G)$-space~$X$ such that~$G$ acts freely and properly.
There, the category of locally compact groupoids is obtained by taking isomorphisms classes of correspondences as
arrows; isomorphic objects in this category are called {\it equivalent}.
We shall see that the f\/inal def\/inition of a~correspondence which will be adopted in this paper is very close to this
initial def\/inition: we shall def\/ine a~correspondence from~$H$ to~$G$ as an~$(H,G)$-space~$X$ such that~$G$ acts properly
(we shall also require some extra data).
For consistency throughout the paper, left and right have been exchanged with respect to some original def\/initions.
The reader should keep in mind that a~correspondence from~$H$ to~$G$ should def\/ine a~$C^*$-correspondence from $C^*(H)$
to $C^*(G)$ and induce representations from~$G$ to~$H$.

The notion of a~Hilsum--Skandalis map, introduced in~\cite{hs:morphisms} and emphasized by J.~Mr\v{c}un
in~\cite{mrc:1999} is akin to the above def\/inition.
A~{\it principal bibundle} is an~$(H,G)$-space~$X$ such that~$H$ acts freely and properly and $s_X:X\rightarrow G^{(0)}$
identif\/ies $H\backslash X$ and $G^{(0)}$.
This is also called a~gene\-ra\-li\-zed homomorphism since a~groupoid homomorphism from~$G$ to~$H$ def\/ines an~$(H,G)$-principal
bibundle.
A~Hilsum--Skandalis map is def\/ined as an isomorphism class of principal bibundles.
In the case when~$H$ and~$G$ are \'etale and smooth, Mr\v{c}un associates to a~principal $(H,G)$-bibundle
a~$(C_c^\infty(H),C_c^\infty(G))$-bimodule and shows the functoriality of the construction.
He also introduces a~notion of Morita equivalence for these smooth algebras.
These are algebraic constuctions.
Let us turn to the $C^*$-algebraic framework and recall the now classical notion of $C^*$-correspondence.

\begin{Definition}
Let~$A$ and~$B$ be C$^*$-algebras.
An {\it $(A,B)$-$C^*$-correspondence} is a~right~$B$-C$^*$-module $\mathcal E$ together with a~$*$-homomorphism
$\pi:A\rightarrow {\mathcal L}_B({\mathcal E})$.
\end{Definition}

As is well-known (see~\cite[Theorem~5.9]{rie:induced}), $C^*$-correspondences can be composed and strong Morita
equivalences are exactly the invertible $C^*$-correspondences.
The category of $C^*$-algebras is def\/ined as the category whose objects are the $C^*$-algebras and arrows are the
isomorphism classes of $C^*$-correspondences.
A~f\/irst step in the generalization of the Muhly--Renault--Williams theorem~\cite[Theorem~2.8]{mrw:Morita} from
equivalences to correspondences was made by M.~Macho Stadler and M.~O'uchi in~\cite{mm:1999}.
Their def\/inition is equivalent to the following: a~correspondence from~$H$ to~$G$ in the sense of Macho Stadler and
O'uchi is an~$(H,G)$-space~$X$ such that both~$G$ and~$H$ act properly and $s_X:X\rightarrow G^{(0)}$ identif\/ies
$H\backslash X$ and $G^{(0)}$.
Assuming that~$H$ and~$G$ have Haar systems, they construct a~$C^*$-correspondence, which we shall denote by $C_r^*(X)$,
from the reduced $C^*$-algebra $C_r^*(H)$ to the reduced $C^*$-algebra $C_r^*(G)$.
N.~Landsman gives in~\cite[Section~3.3]{lan:2001} the same construction, in the category of Lie groupoids.
In his def\/inition of a~correspondence from~$H$ to~$G$, he assumes that the action of~$H$ is free and proper and the
action of~$G$ is proper.
He constructs the $C^*$-correspondence $C^*(X)$ from $C^*(H)$ to $C^*(G)$ and sketches the proof of the functoriality of
this construction.
In~\cite[Theorem~7.11]{tu:2004}, J.-L.~Tu gives the same result, but assuming only that~$G$ and~$H$ are locally compact,
not necessarily Hausdorf\/f, groupoids with Haar systems.
He states without proof the functoriality of the construction.

The notion of correspondence used by M.~Macho Stadler, M.~O'uchi et al.\ is insuf\/f\/icient.
Indeed, it does not cover the restriction of a~representation of a~locally compact group~$G$ to a~closed subgroup~$H$.
In that case $X=G$ as an~$(H,G)$-space but the assumption that $H\backslash G=G^{(0)}$ does not hold.
Restriction of representations is covered by the following construction.
In~\cite{bs:morphisms}, M.~Buneci and P.~Stachura def\/ine a~morphism from~$H$ to~$G$ as an $(H,G)$-space~$X$, where~$X$
is~$G$ itself and~$G$ acts on~$X$ by right multiplication.
Note that this is a~particular case of our initial def\/inition of a~correspondence.
Assume moreover that~$H$ and~$G$ are locally compact and have Haar systems.
It is not dif\/f\/icult to def\/ine an action of $C_c(H)$ on $C_c(G)$ by multipliers.
However, in order to have a~$*$-homomorphism from $C^*(H)$ into the multiplier algebra $MC^*(G)$, it is necessary to
introduce a~modular function.
Since this condition is fairly technical and we are going to give it in a~more general framework, we do not reproduce it
here.
Buneci and Stachura show also that the construction of the $C^*$-algebra is functorial, where the arrows in the category
of $C^*$-algebras are the $*$-homomorphisms $A\rightarrow MB$.
One of the drawbacks of their construction is that it does not include groupoid equivalences.

\subsection{Holkar's def\/inition}\label{section2.2}

Let us now present the notion of correspondence given by R.D.~Holkar in his Ph.D.~thesis~\cite{hol:thesis}, which includes both the correspondences in the sense of Macho Stadler and~O'uchi and the
morphisms of Buneci and Stachura.
We f\/irst introduce a~notation and a~def\/inition.

Let~$G$ be a~locally compact groupoid and~$X$ a~locally compact left~$G$-space.
We form the semi-direct product groupoid $G\ltimes X$.
A~left Haar system~$\lambda$ for~$G$ def\/ines $\lambda_1$, integration along the f\/ibers of the range map $r:G\ltimes
X\rightarrow X$ given by $r(\gamma, x)=\gamma x$ and $\lambda_2$, integration along the f\/ibers of the source map
$s:G\ltimes X\rightarrow X$, given by $s(\gamma,x)=x$:
\begin{gather*}
\lambda_1(f)(x)=\int f\big(\gamma,\gamma^{-1}x\big)d\lambda^{r(x)}(\gamma),
\qquad
\lambda_2(f)(x)=\int f\big(\gamma^{-1}, x\big)d\lambda^{r(x)}(\gamma)
\end{gather*}
for $f\in C_c(G\ltimes X)$.

\begin{Definition}
\label{Delta}
Let $\Delta:G\ltimes X\rightarrow{\bf{R}}_+^*$ be a~Borel cocycle.
We say that a~measure~$\mu$ on~$X$ is a~$\Delta$-measure with respect to $(G,\lambda)$ if
\begin{gather*}
\mu\circ\lambda_1=\Delta(\mu\circ\lambda_2).
\end{gather*}
In other words,~$\mu$ is a~quasi-invariant measure and admits~$\Delta$ as Radon--Nikodym derivative.
\end{Definition}

\begin{Definition}[\cite{hol:thesis}]
\label{correspondence}
Let $(G,\lambda)$ and $(H,\beta)$ be locally compact groupoids with Haar systems.
We say that $(X,\alpha)$ is a~$\Delta$-correspondence from $(H,\beta)$ to $(G,\lambda)$ if
\begin{enumerate}\itemsep=0pt
\item $X$ is a~locally compact space $(H,G)$-space;
\item the right action of~$G$ is proper;
\item $\Delta:H\ltimes (X/G)\rightarrow{\bf{R}}_+^*$ is a~continuous cocycle;
\item $\alpha=(\alpha_u)_{u\in G^{(0)}}$ is a~continuous~$G$-equivariant system of measures for $s_X:X\rightarrow
G^{(0)}$;
\item for all $u\in G^{(0)}$, $\alpha_u$ is a~$\Delta$-measure with respect to $(H,\beta)$.
\end{enumerate}
We say that~$\Delta$ is the module of the correspondence $(X,\alpha)$.
\end{Definition}

\subsection{Examples}\label{examples}

\begin{enumerate}\itemsep=0pt
\item[a)] In the case of a~correspondence in the sense of Macho Stadler and O'uchi, one def\/ines the $s_X$-system~$\alpha$ by
\begin{gather*}
\alpha(f)(u)=\int f\big(h^{-1}x\big)d\beta^{r_X(x)}(h),
\end{gather*}
where $f\in C_c(X)$, $u\in G^{(0)}$ and one has chosen $x\in X$ such that $s_X(x)=u$; the integral depends only on~$u$.
In that case $\Delta\equiv1$.
Note that this also includes groupoid equivalences, which are a~particular case of Macho Stadler and O'uchi
correspondences.

\item[b)] In the case of a~morphism in the sense of Buneci and Stachura where $X=G$, one def\/ines the $s_X$-system~$\alpha$ by
\begin{gather*}
\alpha(f)(u)=\int f\big(\gamma^{-1}\big)d\lambda^u(\gamma),
\qquad
f\in C_c(G).
\end{gather*}
Their condition is exactly the assumption~(v) of Def\/inition~\ref{correspondence}.

\item[c)] Let~$H$ be a~closed subgroup of a~locally compact group~$G$ endowed respectively with left Haar measures~$\beta$
and~$\lambda$.
The correspondence which gives induction of a~representation of~$H$ to a~representation of~$G$ is $X=G$ as
a~$(G,H)$-space, endowed with the right Haar measure $\alpha=\lambda^{-1}$.
An easy computation shows that~$\alpha$ is a~$1$-measure with respect to $(G,\lambda)$.
On the other hand the correspondence which gives the restriction of a~representation of~$G$ to a~representation of~$H$
is $Y=G$ as a~$(H,G)$-space, again endowed with the right Haar measure $\alpha=\lambda^{-1}$.
Another easy computation shows that~$\alpha$ is a~$\Delta$-measure with respect to $(H,\beta)$, where
$\Delta(h)=\delta_H(h)/\delta_G(h)$ for $h\in H$ and $\delta_H$, $\delta_G$ are respectively the modular functions of~$H$
and~$G$.
\end{enumerate}

\subsection[The $C^*$-correspondence]{The $\boldsymbol{C^*}$-correspondence}\label{section2.4}

The construction of the $C^*$-correspondence $C^*(X,\alpha)$ from a~groupoid $(H,G)$-correspondence $(X,\alpha)$ as
above is straightforward.
We refer to~\cite{hol:thesis} (or to the previous constructions) for details.
One f\/irst ignores~$H$ and constructs the right $C^*$-module ${\mathcal E}=C^*(X,\alpha)$ over $B=C^*(G,\lambda)$ as
in~\cite{ren:representations} or~\cite{hr:2013}; it is the completion of $C_c(X)$ for the following operations: for
$f\in C_c(G)$ and $\xi,\eta\in C_c(X)$,
\begin{gather*}
\xi f(x)=\int \xi(x\gamma)f\big(\gamma^{-1}\big)d\lambda^{s(x)}(\gamma),
\qquad
\langle\xi,\eta\rangle(\gamma)=\int \overline{\xi(x)} \eta(x\gamma) d\alpha_{r(\gamma)}(x).
\end{gather*}
The $*$-homomorphism $\pi:C^*(H,\beta)\rightarrow {\mathcal L}_B({\mathcal E})$ is def\/ined, for $g\in C_c(H)$ and
$\xi\in C_c(X)$ by
\begin{gather*}
g\xi(x)=\int g(h)\xi\big(h^{-1}x\big)\Delta^{1/2}\big(h^{-1}, x\big)d\beta^{r(x)}(h).
\end{gather*}
\begin{Lemma}
For $g\in C_c(H)$ and $\xi,\eta\in C_c(X)$, one has
$\langle g\xi,\eta\rangle=\langle\xi, g^*\eta \rangle$.
\end{Lemma}
\begin{proof}
This is a~simple computation which justif\/ies the introduction of the module~$\Delta$.
\end{proof}

By def\/inition of the full $C^*$-norm,~$\pi$ extends to $C^*(H,\beta)$.
This gives the $C^*$-correspondence $C^*(X,\alpha)$ from $C^*(H,\beta)$ to $C^*(G,\lambda)$.

\subsection{Composition of correspondences}\label{section2.5}

We sketch here the construction of the composition of correspondences.
A~detailed presentation is given in~\cite{hol:thesis}.
The construction leans upon two elementary results which we f\/irst recall.

\begin{Lemma}
\label{trivialization}
Let $c:G\rightarrow{\bf{R}}$ be a~continuous cocycle, where~$G$ is a~locally compact groupoid which has a~Haar system.
Assume that the groupoid~$G$ is proper.
Then, there exists a~continuous function $b: G^{(0)}\rightarrow{\bf{R}}$ such that $c(\gamma)=b\circ r(\gamma)-b\circ s(\gamma)$.
\end{Lemma}
\begin{proof}
Let $\pi: G^{(0)}\rightarrow G\backslash G^{(0)}$ be the quotient map.
Since $G\backslash G^{(0)}$ is paracompact and~$\pi$ is an open map, we can apply~\cite[Appendice~1, Lemme~1]{bou:int7-8}:
there exists $F: G^{(0)}\rightarrow {\bf{R}}^+$ continuous such that:
\begin{enumerate}\itemsep=0pt
\item[a)] $F$ is not identically zero on any equivalence class modulo~$G$;
\item[b)] for every compact set $K\subset G\backslash G^{(0)}$, the intersection of $\pi^{-1}(K)$ with the support
of~$F$ is compact.
\end{enumerate}

Let $\lambda=(\lambda^x)$ be a~Haar system for~$G$.
The function $h:G^{(0)}\rightarrow{\bf{R}}^+$ def\/ined by $\displaystyle
h(x)=\int (F\circ s)(\gamma)d\lambda^x(\gamma)$
for $x\in G^{(0)}$ is well def\/ined and strictly positive.
Moreover, it is constant on the equivalence classes.
If we replace~$F$ by $F/h$, we obtain a~function (still denoted by~$F$) which satisf\/ies a), b) and the following
condition:
\begin{enumerate}\itemsep=0pt
\item[c)] $\displaystyle
\int (F\circ s)(\gamma)d\lambda^x(\gamma)=1$ for all $x\in G^{(0)}$.
\end{enumerate}

Such a~function is also constructed in~\cite[Section~6]{tu:2004} where it is called a~cutof\/f function.
Then, the integral
\begin{gather*}
b(x)=\int (F\circ s)(\gamma) c(\gamma)d\lambda^x(\gamma)
\end{gather*}
is well def\/ined and def\/ines a~continuous function such that $c(\gamma)=b\circ r(\gamma)-b\circ s(\gamma)$.
\end{proof}

The technique of the second result is also adapted from~\cite{bou:int7-8}.
It can also be found in~\cite[Appendix~A.1]{adr:amenable} in the case of a~proper Borel equivalence relation.

\begin{Lemma}\label{proper}
Let~$G$ be a~proper locally compact groupoid with Haar system~$\lambda$.
Then,
\begin{enumerate}\itemsep=0pt
\item[$(i)$] the quotient map $\pi: G^{(0)}\rightarrow G\backslash G^{(0)}$ carries a~continuous~$\pi$-system of measures
$\dot\lambda$ defined~by
\begin{gather*}
\dot\lambda(f)(\dot x)=\int (f\circ s)d\lambda^x
\end{gather*}
where $ f\in C_c(G^{(0)})$, $x\in G^{(0)}$ and $\dot x=\pi(x)$;
\item[$(ii)$] a~Radon measure~$\mu$ on $G^{(0)}$ is of the form $m\circ\dot\lambda$ for some Radon measure~$m$ on $G\backslash
G^{(0)}$ if and only if the measure $\mu\circ\lambda$ on~$G$ is symmetric, i.e.~satisfies
$\mu\circ\lambda=(\mu\circ\lambda)^{-1}$.
\end{enumerate}
\end{Lemma}
\begin{proof}
Since this is classical, we just give a~sketch of the proof.
The integral in (i) is well def\/ined because~$G$ is proper.
The necessity of the symmetry of $m\circ\dot\lambda\circ\lambda$ is given by Fubini's theorem.
The construction of~$m$ uses the cutof\/f function~$F$ of the previous lemma.
Explicitly, for $f\in C_c(G\backslash G^{(0)})$, $m(f)=\mu(F(f\circ\pi))$.
\end{proof}
One can combine this lemma and Lemma~\ref{measurable trivialization} given below to describe all quasi-invariant
measures of a~proper groupoid: up to equivalence, they are given by measures on the quotient space.

We can now construct the composition of two correspondences.
Let $(G_i,\lambda_i)$, where $i=1,2,3$, be locally compact groupoids with Haar systems, $(X,\alpha)$ a~correspondence
from $G_1$ to $G_2$ with module $\Delta_X: G_1\ltimes X/G_2\rightarrow{\bf{R}}_+^*$ and $(Y,\beta)$ a~correspondence from
$G_2$ to $G_3$ with module $\Delta_Y: G_2\ltimes Y/G_3\rightarrow{\bf{R}}_+^*$.
One constructs a~correspondence $(Z, \tau)$ from $G_1$ to $G_3$ with module $\Delta: G_1\ltimes
Z/G_3\rightarrow{\bf{R}}_+^*$ as follows.
As a~$(G_1, G_3)$-space,~$Z$ is the usual composition $(X*Y)/G_2$, where $X*Y$ is the f\/ibered product over $G_2^{(0)}$
and $G_2$ acts by $(x,y)h=(xh, h^{-1}y)$.
One def\/ines the system of measures $(\beta_u\circ\alpha)_{u\in G_3^{(0)}}$ for the map $s:X*Y\rightarrow G_3^{(0)}$
sending $(x,y)$ to $s_Y(y)$ by
\begin{gather*}
\int f d(\beta_u\circ\alpha)=\int f(x,y)d\alpha_{r(y)}(x) d\beta_u(y)
\end{gather*}
for $f\in C_c(X*Y)$.
Let $\pi: X*Y\rightarrow Z=X*Y/G_2$ denote the quotient map.
It carries the system of measures $\lambda=(\lambda_z)_{z\in Z}$ def\/ined by
\begin{gather*}
\int f d\lambda_{\pi(x,y)}=\int f\big(xh, h^{-1}y\big)d\lambda^{s(x)}_2(h).
\end{gather*}
Since $G_2$ acts properly on~$X$, the groupoid $(X*Y)\rtimes G_2$ is proper.
We view $\Delta_Y$ as a~cocycle on $(X*Y)\rtimes G_2$ and use Lemma~\ref{trivialization} to trivialize it: there exists
a~continuous $b: X*Y\rightarrow{\bf{R}}_+^*$ such that for all $(x,y)\in X*Y$ and $h\in G_2$ such that $s(h)=r(y)$, the
equality
\begin{gather*}
\Delta_Y(h,y)=b(xh^{-1}, hy)/b(x,y)
\end{gather*}
holds.
Moreover, since $\Delta_Y(h,y)$ depends only on the class $\dot y$ in $Y/G_3$, we may assume that $b(x,y)$ depends only
on $(x,\dot y)$.
Let us f\/ix $u\in G_3^{(0)}$.
The measure $b(\beta_u\circ\alpha)$ on $X*Y$ satisf\/ies the condition (ii) of Lemma~\ref{proper} with respect to
$((X*Y)\rtimes G_2,\lambda_2)$.
Hence there exists a~measure $\tau_u$ on $Z=X*Y/G_2$ such that
\begin{gather*}
b(\beta_u\circ\alpha)=\tau_u\circ \lambda.
\end{gather*}
Note that $\tau_u$ is supported on $Z_u=\pi(X*s_Y^{-1}(u))$.
This def\/ines the desired system $\tau=(\tau_u)_{u\in G_3^{(0)}}$.
Since the measures $b(\beta_u\circ\alpha)$ are~$\Delta$-measures with respect to $(G_1,\lambda_1)$, where
\begin{gather*}
\Delta(\gamma, (x,y))=b(\gamma x,y)^{-1}\Delta_X(\gamma,x) b(x,y),
\end{gather*}
one deduces (see Proposition~\ref{key} below) that $\Delta(\gamma, (x,y))$ depends only on $(\gamma, \pi (x,y))$ and
that the measures $\tau_u$ are~$\Delta$-measures with respect to $(G_1,\lambda_1)$.
We write $Z=X*_{G_2}Y$ and $\tau=\alpha*_{G_2}\beta$.

We state without a~proof the expected result:

\begin{Theorem}[\cite{hol:thesis}]
The construction of the $C^*$-algebra of a~locally compact groupoid with Haar system is functorial.
More precisely, given $(G_i,\lambda_i)$, where $i=1,2,3$, $(X,\alpha)$ and $(Y,\beta)$ as above, the
$C^*$-correspondence $C^*(X*_{G_2}Y,\alpha*_{G_2}\beta)$ is isomorphic to the composition of the $C^*$-correspondences
$C^*(X,\alpha)\otimes_{C^*(G_2,\lambda_2)} C^*(Y,\beta)$.
\end{Theorem}

\section{Induced representations}\label{section3}

As we have seen, a~correspondence $(X,\alpha)$ from $(H,\beta)$ to $(G,\lambda)$ gives a~$C^*$-correspondence
$C^*(X,\alpha)$ from $C^*(H,\beta)$ to $C^*(G,\lambda)$, hence it induces a~map from ${\rm Rep}(G)$ to ${\rm Rep}(H)$,
where ${\rm Rep}(G)$ is the set of equivalence classes of representations of~$G$.
The purpose of this section is to describe this map without passing through the $C^*$-algebras.
Recall from~\cite{ren:representations} that a~representation of~$G$ is a~pair $(m, {\mathcal H})$, where~$m$ is
a~transverse measure class on~$G$ and ${\mathcal H}$ is a~measurable~$G$-Hilbert bundle.
Since we have f\/ixed a~Haar system~$\lambda$ for~$G$, the transverse measure class~$m$ is given by a~quasi-invariant
measure~$\mu$ on $G^{(0)}$.
We denote by $\Delta_\mu: G\rightarrow {\bf{R}}_+^*$ its module.

\subsection{Disintegration of quasi-invariant measures}\label{section3.1}

The following proposition is a~slight variation of~\cite[Corollary 5.3.11]{adr:amenable}.
It generalizes~\cite[Subsection~I.3.21]{ren:approach} and plays a~crucial role in the construction of induced representations.
We use the notation $\lambda_1$, $\lambda_2$ introduced before Def\/inition~\ref{Delta}.

\begin{Proposition}
\label{key}
Let $(G,\lambda)$ be a~locally compact groupoid with Haar system.
Let $X$, $Y$ be left $G$-spaces and $\pi: X\rightarrow Y$ be a~continuous~$G$-equivariant surjection.
Let~$\nu$ be a~$\sigma$-finite measure on~$X$ and let $\displaystyle
\nu=\int \rho_y d\mu(y)$ a~disintegration of~$\nu$ along~$\pi$.
\begin{enumerate}\itemsep=0pt
\item[$(i)$] If~$\nu$ is quasi-invariant with respect to $(G,\lambda)$, then~$\mu$ is quasi-invariant with respect to
$(G,\lambda)$ and for $(\mu\circ\lambda_1^Y)$-a.e.
$(\gamma, y)\in G*Y$, $\gamma \rho_y\sim \rho_{\gamma y}$.
More precisely, let $\Delta_X: G\ltimes X\rightarrow{\bf{R}}_+^*$ be such that
\begin{gather}
\nu\circ\lambda_1^X=\Delta_X \big(\nu\circ\lambda_2^X\big).
\label{eq:1}
\end{gather}
Then there exist measurable functions $\Delta_Y: G\ltimes Y\rightarrow{\bf{R}}_+^*$ and $\delta: G\ltimes
X\rightarrow{\bf{R}}_+^*$ such that
\begin{gather}
\mu\circ\lambda_1^Y=\Delta_Y \big(\mu\circ\lambda_2^Y\big)
\label{eq:2},
\\
\gamma \rho_y=\delta \big(\gamma^{-1},\cdot\big)\rho_{\gamma y}
\label{eq:3}
\end{gather}
and they are related by
\begin{gather}
\Delta_X(\gamma, x)=\delta(\gamma, x)\Delta_Y(\gamma,\pi(x)).
\label{eq:4}
\end{gather}
for $\nu\circ\lambda_1^X$ a.e.
$(\gamma,x)\in G*X$.
\item[$(ii)$] Conversely, suppose that there exist measurable functions $\Delta_Y: G\ltimes Y\rightarrow{\bf{R}}_+^*$ and $\delta:
G\ltimes X\rightarrow{\bf{R}}_+^*$ such that~\eqref{eq:2} and~\eqref{eq:3} hold.
Then $\nu=\mu\circ\rho$ is a~$\Delta_X$-measure, where $\Delta_X$ is defined by~\eqref{eq:4}.
\end{enumerate}
\end{Proposition}

\begin{proof}
We start with an observation.
Let us introduce the map
\begin{gather*}
{\rm id}\times\pi: \ G\ltimes X\rightarrow G\ltimes Y.
\end{gather*}
We endow it with the system of measures $1\otimes \rho$.
By construction, we have the commutation relation
\begin{gather*}
\lambda_2^Y\circ (1\otimes\rho)=\rho\circ\lambda_2^X.
\end{gather*}
However, a~similar commutation relation for $\lambda_1^X$ and $\lambda_1^Y$ requires the relation~\eqref{eq:3}.
More precisely, if~\eqref{eq:3} holds, a~lengthy but straightforward computation gives
\begin{gather*}
\lambda_1^Y \circ(1\otimes\rho)=\delta^{-1}\big(\rho\circ\lambda_1^X\big).
\end{gather*}
Let us f\/irst prove the assertion (ii).
We assume that~\eqref{eq:2} and~\eqref{eq:3} hold.
Starting from the equality
\begin{gather*}
\nu\circ\lambda^X_1=\mu\circ\rho\circ\lambda^X_1,
\end{gather*}
and identifying functions and the corresponding multiplication operators, we have
\begin{gather*}
\nu\circ\lambda^X_1\circ\delta^{-1}=\mu\circ\rho\circ\lambda^X_1\circ\delta^{-1}
=\mu\circ\lambda_1^Y \circ(1\otimes\rho)
=\mu\circ\lambda_2^Y\circ \Delta_Y \circ(1\otimes\rho)
\\
\phantom{\nu\circ\lambda^X_1\circ\delta^{-1}}
=\mu\circ\lambda_2^Y\circ(1\otimes\rho)\circ (\Delta_Y \circ({\rm id}\times\pi))
=\mu\circ \rho\circ\lambda_2^X\circ (\Delta_Y \circ({\rm id}\times\pi))
\\
\phantom{\nu\circ\lambda^X_1\circ\delta^{-1}}
=\nu\circ\lambda_2^X\circ (\Delta_Y \circ({\rm id}\times\pi)).
\end{gather*}
Thus, we can write $\nu\circ\lambda_1^X=\Delta_X(\nu\circ\lambda_2^X)$ where
\begin{gather*}
\Delta_X\delta^{-1}=\Delta_Y \circ({\rm id}\times\pi)
\qquad
\text{or equivalently}
\qquad
\Delta_X(\gamma,x)=\Delta_Y(\gamma,\pi(x))\delta(\gamma,x).
\end{gather*}
This shows that~$\nu$ is a~$\Delta_X$-measure with respect to $(G,\lambda)$.

Let us prove the assertion (i).
We assume that~\eqref{eq:1} holds.
As pseudo-images of the equivalent measures $\mu\circ\lambda_2^X$ and $\mu\circ\lambda_1^X$ under $id\times\pi: G\ltimes
X\rightarrow G\ltimes Y$, the measures $\mu\circ\lambda_2^Y$ and $\mu\circ\lambda_1^Y$ are equivalent.
Therefore there exists a~measurable function $\Delta_Y: G\ltimes Y\rightarrow {\bf{R}}_+^*$ such that~\eqref{eq:2} holds.
We are going to compare two disintegrations of the measure $\nu\circ\lambda_1^X$ on $G\ltimes X$ along $id\times\pi: G\ltimes
X\rightarrow G\ltimes Y$, taking the same measure $\mu\circ\lambda_1^Y$ as base measure.
The f\/irst one is obtained by applying the inverse map of $G\ltimes X$ to the disintegration
\begin{gather*}
\nu\circ\lambda_2^X=\big(\mu\circ\lambda_2^Y\big)\circ (1\otimes\rho).
\end{gather*}
This yields after some computation:
\begin{gather*}
\nu\circ\lambda_1^X=\big(\mu\circ\lambda_1^Y\big)\circ \Phi
\qquad
\text{where}
\qquad
\Phi(f)(\gamma, y)=\int f\big(\gamma,\gamma^{-1}x\big)d\rho_{\gamma y}(x).
\end{gather*}
The second disintegration is obtained from the relations~\eqref{eq:1} and~\eqref{eq:2}:
\begin{gather*}
\nu\circ\lambda_1^X=\nu\circ\lambda_2^X\circ\Delta_X
=\mu\circ\rho\circ\lambda_2^X\circ\Delta_X
=\mu\circ\lambda_2^Y\circ(1\otimes\rho)\circ\Delta_X
\\
\phantom{\nu\circ\lambda_1^X}
=\mu\circ\lambda_1^Y\circ\Delta_Y^{-1}\circ(1\otimes\rho)\circ\Delta_X.
\end{gather*}
By uniqueness of the disintegration along the map $id\times\pi$, we get
\begin{gather*}
\Phi=\Delta_Y^{-1}\circ(1\otimes\rho)\circ\Delta_X
\end{gather*}
or equivalently, for all $f\in C_c(G\ltimes X)$ and $\nu\circ\lambda_1^Y$ a.e.
$(\gamma,y)\in G\ltimes Y$,
\begin{gather*}
\int f\big(\gamma,\gamma^{-1}x\big)d\rho_{\gamma y}(x)=\int f(\gamma,x){\Delta_X(\gamma,x)\over
\Delta_Y(\gamma,\pi(x))}d\rho_y(x).
\end{gather*}
This gives the formula
\begin{gather*}
\gamma^{-1}\rho_{\gamma y}=\delta(\gamma,\cdot)\rho_y,
\qquad \text{where} \qquad  \delta(\gamma,x)=\frac{\Delta_X(\gamma,x)}{\Delta_Y(\gamma,\pi(x))}.\tag*{\qed}
\end{gather*}
  \renewcommand{\qed}{}
\end{proof}

\subsection{Construction of the induced quasi-invariant measure}\label{section3.2}

Let us f\/irst state and prove a~result which is the measurable version of Lemma~\ref{trivialization}.

\begin{Lemma}
\label{measurable trivialization}
Let $c:G\rightarrow{\bf{R}}$ be a~measurable cocycle, where~$G$ is a~locally compact groupoid endowed with a~Haar system
and a~quasi-invariant measure.
Assume that the groupoid~$G$ is proper.
Then, there exists a~measurable function $b: G^{(0)}\rightarrow{\bf{R}}$ such that $c(\gamma)=b\circ r(\gamma)-b\circ
s(\gamma)$ for all $\gamma\in G$.
\end{Lemma}
\begin{proof}
We denote by~$R$ be the graph of the equivalence relation on $G^{(0)}$ def\/ined by~$G$ and by $G'$ the isotropy bundle
of~$G$, i.e.~the subgroupoid def\/ined by $r(\gamma)=s(\gamma)$ (it is the union of the isotropy groups).
Since~$G$ is proper, the isotropy groups $G(x)$, where $x\in G^{(0)}$, are compact.
Therefore, the restriction of~$c$ to $G(x)$ and to $G'$ is trivial.
Since~$R$ can be identif\/ied with the quotient $G'\backslash G$, there exists a~measurable cocycle ${\underline c}:
R\rightarrow{\bf{R}}$ such that for all $\gamma\in G$, $c(\gamma)={\underline c}(r(\gamma), s(\gamma))$.
Since~$R$ is a~closed subset of $G^{(0)}$, the quotient map $\pi: G^{(0)}\rightarrow G\backslash G^{(0)}$ has a~Borel
section~$\sigma$ (see~\cite[Theorem~2.1]{ram:dichotomy}).
The function $b: G^{(0)}\rightarrow{\bf{R}}$ def\/ined by $b(x)={\underline c}(x,\sigma(x))$ is a~measurable coboundary
for~$c$.
\end{proof}

We now begin the construction of the representation $(\underline\mu, \underline{\mathcal H})$ of $(H,\beta)$ induced by
the representation $(\mu, {\mathcal H})$ of $(G,\lambda)$ through the correspondence $(X,\alpha)$.
We f\/irst def\/ine the quasi-invariant measure $\underline\mu$.
The module of~$\mu$ is denoted by $\Delta_\mu$.
Using~\cite[Theorem~3.2]{ram:topologies}, we choose $\Delta_\mu$ so that it is a~strict homomorphism.
According to Proposition~\ref{key} applied to the~$G$-map $s_X:X\rightarrow G^{(0)}$, the measure $\mu\circ\alpha$
on~$X$ is a~$\Delta_X$-measure on~$X$ with respect to $(G,\lambda)$, where $\Delta_X(x,\gamma)=\Delta_\mu(\gamma)$.
We use the same idea as in Section~\ref{section2.5}, namely we trivialize the cocycle $\Delta_X$ on the proper groupoid $X\rtimes G$.
We apply Lemma~\ref{measurable trivialization} to obtain a~measurable function $b: X\rightarrow {\bf{R}}_+^*$ such that
$\Delta_X(x,\gamma)=b(x\gamma)/b(x)$.
Then, the measure $b(\mu\circ\alpha)$ is an invariant measure with respect to $(X\rtimes G, \lambda)$.
According to Lemma~\ref{proper}, there is a~measure~$m$ on $X/G$ such that
\begin{gather*}
b(\mu\circ\alpha)=m\circ\dot\lambda,
\end{gather*}
where we have used the same notation as in the lemma.
Here,
\begin{gather*}
\dot\lambda(f)(\dot x)=\int f(x\gamma)d\lambda^{s(x)}(\gamma)
\qquad
\text{for}
\quad
f\in C_c(X).
\end{gather*}
Next, we study the invariance property of~$m$ with respect to $(H,\beta)$.
Since the $\alpha_u$'s are~$\Delta$-measures, so is $\mu\circ\alpha$.
This implies that the measure $m\circ\dot\lambda$ satisf\/ies:
\begin{gather*}
m\circ\dot\lambda\circ\beta_1={b\circ r\over b\circ s}\Delta(m\circ\dot\lambda\circ\beta_2).
\end{gather*}
Since the system $\dot\lambda$ is invariant under~$H$, Proposition~\ref{key} gives that
\begin{gather*}
m\circ\beta_1=\Delta_m(m\circ\beta_2),
\end{gather*}
where
\begin{gather}\label{eq:6}
\Delta_m(h,\dot x)={b(h x)\over b(x)}\Delta(h,\dot x).
\end{gather}

The last step is to pass from the measure~$m$ on $X/G$ to the measure $\underline\mu$ on $H^{(0)}$.
We just choose a~pseudo-image $\underline\mu$ of~$m$ and disintegrate~$m$ along $\dot r:X/G\rightarrow H^{(0)}$:
\begin{gather}\label{eq:7}
m=\int\rho^u d\underline\mu(u).
\end{gather}
According to Proposition~\ref{key}, $\underline\mu$ is quasi-invariant and we have cocycles $\Delta_{\underline\mu}:
H\rightarrow{\bf{R}}_+^*$ and $\delta: H\ltimes X/G\rightarrow{\bf{R}}_+^*$ such that $\underline\mu$ is quasi-invariant
with module $\Delta_{\underline\mu}$ and
\begin{gather}\label{eq:8}
h\rho^{s(h)}=\delta(h,\cdot)\rho^{r(h)},
\\
\label{eq:9}
\Delta_{\underline\mu}(h)=\Delta_m(h, \dot x)\delta(h,\dot x).
\end{gather}

\subsection{Construction of the induced Hilbert bundle}\label{section3.3}

Given a~representation $(\mu, {\mathcal H})$ of $(G,\lambda)$, we have constructed in the previous subsection the
induced quasi-invariant measure $\underline\mu$ of $(H,\beta)$.
Let us now construct the induced~$H$-Hilbert bundle $\underline{\mathcal H}$.
For $x\in X$ we denote by $G(x)=\{\gamma\in G: \gamma x=x\}$ the stabilizer of~$x$ and by $\kappa_x$ its normalized Haar
measure.
Then $\kappa=(\kappa_x)_{x\in X}$ is a~Borel Haar system for the the isotropy group bundle $G(X)$ of $X\rtimes G$ (see for
example~\cite[Lemma~1.5]{ren:ideal}).
We view the measurable~$G$-Hilbert bundle $X*{\mathcal H}$ over~$X$ as a~$X\rtimes G$-Hilbert bundle.

\begin{Definition}
Let $\mathcal H$ be a~measurable~$G$-Hilbert bundle with bundle map $p:{\mathcal H}\rightarrow G^{(0)}$.
We def\/ine its f\/ixed-point bundle as
\begin{gather*}
{\mathcal H}^G=\{\xi\in {\mathcal H}: \gamma p(\xi)=p(\xi)\Rightarrow L(\gamma)\xi=\xi\}.
\end{gather*}
\end{Definition}
We shall use the following easy result.

\begin{Proposition}
Let $\mathcal H$ be a~measurable~$G$-Hilbert bundle.
\begin{itemize}\itemsep=0pt
\item ${\mathcal H}^G$ is a~measurable~$G$-Hilbert subbundle of $\mathcal H$.
\item If~$G$ is a~proper groupoid, integration over the isotropy subgroups defines a~bundle projection of ${\mathcal H}$
onto ${\mathcal H}^G$.
\end{itemize}
\end{Proposition}

We def\/ine integration over the isotropy subgroup $G(x)$, where $x\in G^{(0)}$, as the orthogonal projection of
${\mathcal H}_x$ onto ${\mathcal H}_x^{G(x)}$
\begin{gather*}
P_x\xi=\int L(\gamma)\xi d\kappa_x(\gamma)
\qquad
\text{if}
\quad
\xi\in {\mathcal H}_x,
\end{gather*}
where $\kappa_x$ is the normalized Haar measure of $G(x)$.
The family $P=(P_x)_{x\in G^{(0)}}$ def\/ines the bundle projection.

Going back to our previous situation, we consider the f\/ixed-point bundle ${\mathcal H}'=(X*{\mathcal H})^{X\rtimes G}$ of
$X*{\mathcal H}$.
The f\/iber above $x\in X$ is:
\begin{gather*}
{\mathcal H}'_x=\{\xi\in {\mathcal H}_{s(x)}: \forall\, h\in G(x), \, L(h)\xi=\xi\}.
\end{gather*}
This bundle is stable under the action of~$G$ and is a~proper~$G$-space.
We denote by ${\mathcal K}={\mathcal H}'/G$ the quotient space.

\begin{Proposition}
The quotient space ${\mathcal K}={\mathcal H}'/G$ is a~measurable Hilbert bundle over $X/G$, where the bundle map
$p:{\mathcal K}\rightarrow X/G$ sends the class $[x,\xi]$ of $(x,\xi)$ to the class $\dot x$ of~$x$.
\end{Proposition}

\begin{proof}
A~choice of $x\in X$ def\/ines a~bijection $\varphi_x:{\mathcal H}'_x\rightarrow {\mathcal K}_{\dot x}$ such that
$\varphi_x(\xi)=[x,\xi]$.
We use this bijection to carry to ${\mathcal K}_{\dot x}$ the Hilbert space structure of ${\mathcal H}'_x$.
This Hilbert space structure does not depend on the choice of~$x$ in its class: suppose that $y=x\gamma$.
The isometry $L(\gamma):{\mathcal H}_{s(\gamma)}\rightarrow {\mathcal H}_{r(\gamma)}$ sends ${\mathcal H}'_y$ onto
${\mathcal H}'_x$ and we have $\varphi_y=\varphi_x\circ L(\gamma)$.
Let $P_x:{\mathcal H}_{s(x)}\rightarrow {\mathcal H}'_x$ be the orthogonal projection.
It satisf\/ies $L(\gamma)\circ P_y=P_x\circ L(\gamma)$.
As it is a~subbundle of $X*{\mathcal H}$, ${\mathcal H}'$ has a~measurable Hilbert bundle structure and so has the
quotient bundle ${\mathcal K}$: we shall say that a~section $\xi$ of $p:{\mathcal K}\rightarrow X/ G$ is measurable if
there is a~measurable section~$\sigma$ for the quotient map $X\rightarrow X/ G$ and a~measurable section~$\eta$ of the
bundle ${\mathcal H}'$ such that $\xi(\dot x)=[\eta\circ\sigma(\dot x)]$.
\end{proof}

The Hilbert bundle $X*{\mathcal H}$ is endowed with the unitary~$H$-action: $h (x,\xi)=(h x,\xi)$.
Since this~$H$-action commutes with the action of~$G$, the subbundle ${\mathcal H}'$ is invariant under~$H$ and the
quotient~${\mathcal K}$ is also a~$H$-Hilbert bundle.
The last step is to go from the~$H$-Hilbert bundle $\mathcal K$ over $X/G$ to a~$H$-Hilbert bundle $\underline{\mathcal
H}$ over $H^{(0)}$.
We use the disintegration $\displaystyle \nu=\int\rho^u d\underline\mu(u)$ of the measure~$\nu$ given earlier.
We def\/ine for $u\in H^{(0)}$:
\begin{gather*}
\underline{\mathcal H}_u=L^2(X/ G,\rho^u, {\mathcal K}).
\end{gather*}
As explained earlier, the measures~$\nu$ on $X/G$ and $\underline\mu$ on $H^{(0)}$ are quasi-invariant and there exists
a~measurable function $\delta: H\ltimes X/G\rightarrow{\bf{R}}_+^*$ such that $h \rho^{s(h)}=\delta(h^{-1},\cdot)\rho^{r(h)}$.
The measurable structure of the bundle $(\underline{\mathcal H}_u)_{u\in H^{(0)}}$ is provided by the family of
measurable sections $u\mapsto \xi_u$, where $\xi$ is a~measurable section of the bundle ${\mathcal K}\rightarrow X/G$
and $\xi_u(\dot x)=\xi(\dot x)$ if $r(\dot x)=u$.
For $h\in H$, we def\/ine $\displaystyle L(h): \underline{\mathcal H}_{s(h)}\rightarrow \underline{\mathcal H}_{r(h)}$ by
\begin{gather*}
(L(h)\xi_{s(h)})(\dot x)=\delta^{1/2}\big(h^{-1},\dot x\big)\xi_{s(h)}\big(h^{-1}\dot x\big).
\end{gather*}
These are unitary operators and $\underline{\mathcal H}$ is the desired~$H$-Hilbert bundle.

\subsection{Examples}\label{section3.4}

Let us spell out the above construction in the classical case~\ref{examples}.c, where~$H$ is a~closed subgroup of~$G$
and $X=G$ is viewed as a~correspondence from~$G$ to~$H$ (note that we considered earlier a~correspondence from~$H$
to~$G$).
To construct the measure~$m$ on $G/H$, we f\/ind $b:G\rightarrow{\bf{R}}_+^*$ such that $b(xh)=b(x)\delta_H(h)$ for
$(x,h)\in G\times H$.
According to Lemma~\ref{proper}(ii), we can write $b\lambda^{-1}=m\circ\dot\beta$ where $\displaystyle\dot\beta(f)(\dot
x)=\int f(xh)d\beta(h)$.
Then~$m$ is quasi-invariant with respect to $(G,\lambda)$.
Since $H^{(0)}$ is reduced to one element~$e$, we have $m=\rho^e$ in equation~\eqref{eq:7}.
The equations~\eqref{eq:6} and~\eqref{eq:8} become respectively $\Delta_m(\gamma,\dot x)=b(\gamma x)/b(x)$ and $\gamma
m=\delta(\gamma,\cdot)m$ where according to equation \eqref{eq:9}, $\delta: G\ltimes (G/H)\rightarrow {\bf{R}}_+^*$ satisf\/ies
$\displaystyle\delta_G(\gamma)=\frac{b(\gamma x)}{b(x)}\delta(\gamma,\dot x)$.
The~$G$-Hilbert space induced by a~$H$-Hilbert space ${\mathcal H}$ is $L^2(G/H, m, (G\times {\mathcal H})/H)$.

\looseness=-1
Let $H_1$, $H_2$ be closed subgroups of a~locally compact group~$G$.
We denote by $\beta_1$, $\beta_2$ and~$\delta_1$,~$\delta_2$ the Haar measures and modules respectively.
Then $(X,\alpha)=(G,\lambda^{-1})$ can also be viewed as a~correspondence from $(H_1,\beta_1)$ to $(H_2,\beta_2)$.
As we have seen earlier,~$\alpha$ is a~$\Delta$-measure with respect to $(H_1,\beta_1)$, where
$\Delta(h)=\delta_1(h)/\delta_G(h)$ for $h\in H_1$.
Given a~$H_2$-Hilbert space ${\mathcal H}$, the induced Hilbert space is just as above $\underline{\mathcal
H}=L^2(G/H_2, m, (G\times {\mathcal H})/H_2)$ and the action of $H_1$ is given by the same formula as before.
If we assume that $H_1$ acts properly on $G/H_2$, this representation can be decomposed over $H_1\backslash G/H_2$.
This is best understood at the level of the correspondence~$X$.
Let $\pi: X\rightarrow H_1\backslash G/H_2$ be the quotient map.
For a~double class $d=H_1xH_2$, we also write $X_d=\pi^{-1}(d)=H_1xH_2$.
We view $X_d$ as a~correspondence from $H_1$ to $H_2$.
It is equipped with a~measure $\alpha_d$ coming from a~disintegration of~$\alpha$ along~$\pi$.
One observes that $X_d/H_2$ and $H_1/(H_1\cap xH_2x^{-1})$ are isomorphic $H_1$-spaces (and their measures match).
Therefore, the $H_1$-Hilbert space $\underline{\mathcal H}_d$ induced through $X_d$ from ${\mathcal H}$ is easily
identif\/ied: it is isomorphic to the $H_1$-Hilbert space induced from the representation ${\mathcal H}^x$ of its subgroup
$H_1\cap xH_2x^{-1}$, given by ${\mathcal H}^x={\mathcal H}$ and $L^x(h)=L(x^{-1}h x)$ for $h\in H_1\cap xH_2x^{-1}$.
Then, the $H_1$-Hilbert space $\underline{\mathcal H}$ is a~direct integral of the $\underline{\mathcal H}_d$'s over
$H_1\backslash G/H_2$.
This is the content of G.~Mackey's subgroup theorem~\cite[Theorem~12.1]{mac:induced}.

\section{Hypergroupoids}\label{section4}

G.~Mackey introduced virtual subgroups (and virtual groups) to the purpose of generalizing the theory of induced
representations.
From the $C^*$-algebraic point of view, induction only requires a~$C^*$-correspondence.
Groupoid correspondences and equivalences suggest to consider more general objects, namely hypergroupoids, and their
$C^*$-algebras, to construct induced representations.
Let us see how hypergroupoid $C^*$-algebras enter into our framework.
When~$X$ is a~$(G,H)$-groupoid equivalence,~$X$ as a~left free and proper~$G$-space determines~$H$ up to isomorphism.
Indeed~$H$ is isomorphic to the groupoid $(X*X)/G$.
Moreover, a~continuous~$G$-equivariant system of measures~$\alpha$ on~$X$ def\/ines a~Haar system on $(X*X)/G$.
Let us assume that~$X$ is a~proper, but not necessarily free,~$G$-space endowed with a~continuous~$G$-equivariant system
of measures~$\alpha$.
It is shown in~\cite{hr:2013} that, although $H=(X*X)/G$ is no longer a~groupoid, the $*$-algebra $C_c(H)$ can be
def\/ined by the same formulas as in the case of a~free and proper~$G$-space.
The representations of~$G$ extend to representations of $C_c(H)$ and there is a~least $C^*$-norm making all these
representations continuous.
We denote by $C_G^*(H)$ the corresponding completion of~$C_c(H)$.
One can also complete $C_c(X)$ into a~$C^*$-correspondence $C^*(X)$ from $C^*(G)$ to $C_G^*(H)$.
Thus one can induce representations from $C_G^*(H)$ to $C^*(G)$.
These are well known constructions in the theory of unitary representations of groups.
The usual situation is a~pair $(G,K)$ where~$K$ is a~compact subgroup of a~locally compact group~$G$.
Then $X=G/K$ is a~proper left~$G$-space which has an invariant measure.
The corresponding hypergroupoid $H=(X*X)/G$ is isomorphic to the double coset hypergroup $K\backslash G/K$.
As we have seen, one can induce representations of~$G$ to representations of the $*$-algebra $C_c(H)$.
The converse problem is more interesting: which representations of $C_c(H)$ induce representations of~$G$? These
problems are studied in the framework of induced representations of hypergroups in~\cite{hkk:hyper} and
in~\cite{her:1992, her:1995}, where the author also uses Rief\/fel's $C^*$-algebraic machinery.

The structure of $H=(X*X)/G$, where~$G$ is a~locally compact groupoid and~$X$ is a~(say left) proper~$G$-space endowed
with a~continuous~$G$-equivariant $r_X$-system~$\alpha$ is that of a~locally compact hypergroupoid with Haar system,
which we are going to def\/ine in this section.
Hypergroupoids of that form will be termed {\it spatial hypergroupoids}.
Hypergroups are hypergroupoids whose unit space is reduced to one point.
Double coset spaces $K\backslash G/K$ as above are spatial hypergroups.
The main dif\/ference between the hypergroups considered here and the classical theory of locally compact hypergroups of
Ch.
Dunkl~\cite{dun:hyper}, R.~Jewett~\cite{jew:hyper} and R.~Spector~\cite{spe:hyper} is that we assume the existence of
a~Haar system, which allows the def\/inition of their $C^*$-algebras just as in the case of groups, while the classical
theory deals with measure algebras.
The def\/inition of a~locally compact hypergroupoid with Haar system which will be proposed in this section is a~tentative
one.
It is inspired by the def\/inition of a~hypergroup given by~\cite{dun:hyper, jew:hyper, spe:hyper} and is rather close in
spirit to~\cite{kpc:hyper}.
Its main virtue is to cover the spatial hypergroupoids as above and the locally compact hypergroups which have a~Haar
measure; moreover, it makes the theory a~development of the case of locally compact groupoids with Haar systems.

The idea of the def\/inition is very simple: we take the usual def\/inition of a~locally compact groupoid~$H$ but where the
product of two composable elements $x$, $y$ is no longer a~third element but a~probability measure $x*y$ with compact
support.
Given a~locally compact space~$X$, $P(X)$~denotes its space of probability measures and given $x\in X$, $\delta_x$
denotes the point mass at~$x$.
For the clarity of the exposition, we present f\/irstly the axioms concerning the product and the involution and secondly
the axioms concerning the Haar system.
However, we are only concerned with locally compact hypergroupoids with Haar systems, that is, hypergroupoids which
satisfy the two sets of axioms.

\begin{Definition}
\label{hypergroupoid}
A~locally compact hypergroupoid is def\/ined as a~pair $(H, H^{(0)})$ of locally compact spaces, range and source maps $r,s:H\rightarrow H^{(0)}$ assumed to be continuous, open and surjective, a~continuous injection $i: H^{(0)}\rightarrow
H$ such that $r\circ i$ and $s\circ i$ are the identity map, a~continuous involution ${\rm inv}:h\mapsto h^*$ of~$H$
such that $r\circ{\rm inv}= s$ and a~product map $m:H^{(2)}\rightarrow P(H)$, where $H^{(2)}$ is the set of composable
pairs, i.e.~pairs $(x,y)\in H\times H$ with $s(x)=r(y)$ (one def\/ines similarly $H^{(3)}$) such that:
\begin{enumerate}\itemsep=0pt
\item the support of $m(x,y)$ is a~compact subset of $H^{r(x)}_{s(y)}$;
\item for all $(x,y,z)\in H^{(3)}$, we have $\displaystyle\int m(x,\cdot)dm(y,z)=\int m(\cdot,z)dm(x,y)$;
\item for all $x\in H$, $m(r(x),x)=m(x,s(x))=\delta_x$;
\item for all $(x,y)\in H^{(2)}$, $m(x,y)^*=m(y^*, x^*)$, where $m(x,y)^*$ is the image of the measure $m(x,y)$ by the
involution;
\item $x=y^*$ if and only if the support of $m(x,y)$ meets $i(H^{(0)})$;
\item for all $f\in C_c(H)$ and $\epsilon>0$, there exists a~neighborhood~$U$ of $i(H^{(0)})$ in~$H$ such that
$|f(x)-f(y^*)|\le \epsilon$ as soon as the support of $m(x,y)$ meets~$U$;
\item for all $x\in H$, the left translation operator $L(x)$ def\/ined by
\begin{gather*}
(L(x)f)(y)=f(x^**y):=\int f dm(x^*,y)
\end{gather*}
sends $C_c(H^{s(x)})$ into $C_c(H^{r(x)})$.
\end{enumerate}
\end{Definition}

\begin{Remark}
When the measures $m(x,y)$ are point masses, one retrieves the product and the inverse map of a~groupoid.
Note however that, as explained below, the continuity of the product is not assumed.
We suspect that adding the second set of axioms concerning Haar systems will force the continuity of the product but we
have not checked this.
Our axioms are modelled after the def\/inition of a~locally compact hypergroup, as given in~\cite{jew:hyper} (where it is
called a~convo) or~\cite[Subsection~1.1.2]{bh:hyper}.
They are tailored to f\/it our main class of examples, namely spatial hypergroupoids, described below in
Theorem~\ref{spatial}.
Simple examples of spatial hypergroupoids, for example the hypergroupoid constructed from the action of
${\bf{Z}}/2{\bf{Z}}$ on ${\bf{R}}$ by the map $x\mapsto -x$, show that the product map is not necessarily continuous, in
the sense that for $f\in C_c(H)$, the map $\displaystyle (x,y)\mapsto f(x*y):=\int f dm(x,y)$ may fail to be continuous on $H^{(2)}$.
Nevertheless, there is a~convolution product which turns $C_c(H)$ into an algebra.
Thus, we drop the continuity of the product but we introduce the axioms which make the construction of the $*$-algebra
$C_c(H)$ and its $C^*$-completions possible.
Axioms (i), (iii), (iv), (v) directly generalize those of~\cite{bh:hyper}.
Axiom (ii) expresses the associativity of the product, when it is extended to bounded Radon measures.
Axiom (vi) does not appear explicitly in the usual axioms for hypergroups but is a~consequence of these axioms
\cite[Lemma~4.3.B]{jew:hyper}.
The main reason to introduce this axiom is that it allows the construction of approximate units in the convolution
algebra $C_c(H)$ just as in~\cite[Proposition~2.1.9]{ren:approach}.
It also implies half of the axiom (v), namely if the support of $m(x,y)$ meets $i(H^{(0)})$, then $x=y^*$.
Axiom (vii) means essentially separate continuity of the product.
\end{Remark}

Let us give the axioms about Haar systems.

\begin{Definition}
\label{Haar}
A~Haar system on a~locally compact hypergroupoid~$H$ is a~system of Radon measures $\lambda=(\lambda^u)_{u\in H^{(0)}}$
for the range map such that
\begin{enumerate}\itemsep=0pt
\item for all $f\in C_c(H)$, $\displaystyle u\in H^{(0)}\mapsto \int fd\lambda^u$ is continuous;
\item for all $f,g\in C_c(H)$ and all $x\in H$,
\begin{gather*}
\int f(x*y)g(y)d\lambda^{s(x)}(y)=\int f(y)g(x^**y)d\lambda^{r(x)}(y);
\end{gather*}
\item for all $f,g\in C_c(H)$, $\displaystyle
x\in H\mapsto \int f(x*y)g(y)d\lambda^{s(x)}(y)$ is continuous with
compact support.
\end{enumerate}
\end{Definition}

\begin{Remark}
Assumption (ii) is called the adjoint property; it is formally stronger than the usual left invariance property
\begin{gather*}
\int f(x*y)d\lambda^{s(x)}(y)=\int f(y)d\lambda^{r(x)}(y).
\end{gather*}
It is shown to be equivalent to the left invariance property in the case of hypergroups.
We expect that this equivalence still holds for our hypergroupoids; for our purpose, which is to def\/ine a~$*$-algebra
structure on $C_c(H)$, we prefer to require the adjoint property.
Assumption~(iii) is also designed to turn $C_c(H)$ into an algebra under convolution.
\end{Remark}

Of course, our def\/inition includes locally compact hypergroups which have a~Haar measure (this includes abelian, compact
or discrete hypergroups, the general case is not settled) but our main motivation is to give a~framework to the spatial
hypergroupoids described in the next theorem.
Spatial hypergroupoids and their $C^*$-algebras appear in~\cite{hr:2013} but without the formal def\/inition of
a~hypergroupoid.

\begin{Theorem}
\label{spatial}
Let~$G$ be a~locally compact groupoid with Haar system and~$X$ a~proper~$G$-space~$X$ endowed with
a~continuous~$G$-equivariant system of measures $\alpha=(\alpha^u)_{u\in G^{(0)}}$.
Then $H=(X*X)/G$ is a~locally compact hypergroupoid with Haar system.
\end{Theorem}

\begin{proof}
We assume that~$X$ is a~left~$G$-space and denote by $r: X\rightarrow G^{(0)}$ its moment map.
Recall that~$G$ acts on $X*X$ by the diagonal action $\gamma(x,y)=(\gamma x, \gamma y)$.
The image of $x\in X$ in $X/G$ is denoted here by $[x]$.
Similarly, the image of $(x,y)\in X*X$ in $(X*X)/G$ is denoted by $[x,y]$.
The unit space of~$H$ is $H^{(0)}=X/G$.
The range and source maps are respectively $r([x,y])=[x]$ and $s([x,y])=[y]$.
The identif\/ication map $i:H^{(0)}\rightarrow H$ is $i([x])=[x,x]$.
The involution map is $[x,y]^*=[y,x]$.
The map $X*X*X\rightarrow H^{(2)}$ sending $(x,y,z)$ to $([x,y],[y,z])$ is surjective and we shall write $[x,y,z]$
instead of $([x,y],[y,z])$ the elements of $H^{(2)}$.
Note that $[x',y',z']=[x,y,z]$ if and only there exist $(\gamma,\zeta)\in G*G(y)$ such that $x'=\gamma x$, $y'=\gamma y$
and $z'=\gamma\zeta z$ where $G(y)=\{\zeta\in G: \zeta y=y\}$ is the isotropy group at~$y$.
We denote by $\beta_y$ the normalized Haar measure of the compact group $G(y)$ and, for $[x,y,z]\in H^{(2)}$, by
$m_{[x,y,z]}$ the measure on~$H$ def\/ined by
\begin{gather*}
\int f dm_{[x,y,z]}=\int f[\zeta x,z] d\beta_y(\zeta).
\end{gather*}
It is then a~tedious but straightforward task to check that all the axioms of Def\/inition~\ref{hypergroupoid} are
satisf\/ied.
Here is a~sketch of the proof of (vi).
One introduces a~proper metric~$d$ on~$X$ def\/ining its topology.
If $f\in C_c(H)$, for a~f\/ixed $y\in X$, the function $x\mapsto f[x,y]$ is uniformly continuous on $X^{r(y)}$.
A~compactness argument gives $\eta>0$ such that for all $(x,z)\in X*X$ such that $d(x,z)<\eta$, the inequality
$|f[x,y]-f[z,y]|<\epsilon$ holds for all $y\in X$ such that $r(y)=r(x)$.
The desired set~$U$ is the image of the subset of $X*X$ def\/ined by $d(x,z)<\eta$.
The system of measures~$\alpha$ def\/ines a~system of measures~$\lambda$ for the range map $H\rightarrow X/G$: for $x\in
X$, $\lambda^{[x]}$ is the image of $\alpha^{r(x)}$ by the proper surjective map $\varphi^x: X^{r(x)}\rightarrow
H^{[x]}$ sending~$y$ to $[x,y]$: for $f\in C_c(H)$,
\begin{gather*}
\int f d\lambda^{[x]}=\int f([x,y])d\alpha^{r(x)}(y).
\end{gather*}
One checks that this depends only on $[x]$.
The continuity of~$\alpha$ ensures the continuity of~$\lambda$.
The disintegration of the measure $\alpha^{r(x)}$ along $\varphi^x: X^{r(x)}\rightarrow H^{[x]}$ is given by
\begin{gather*}
\alpha^{r(x)}=\int \beta^x_{[x,y]}d\lambda^{[x]}([x,y]),
\qquad
\text{where for}
\quad
g\in C_c\big(X^{r(x)}\big)
\quad
\int g d\beta^x_{[x,y]}=\int g(\zeta y)d\beta_x(\zeta).
\end{gather*}
One deduces the expected expression of the convolution product: for $f,g\in C_c(H)$, one has
\begin{gather*}
f*g[x,z]= \int f[x,y]\left(\int g dm_{[y,x,z]}\right)d\lambda^{[x]}([x,y])
\\
\phantom{f*g[x,z]}
= \int f[x,y]\left(\int g[\zeta y,z] d\beta_x(\zeta)\right)d\lambda^{[x]}([x,y])
\\
\phantom{f*g[x,z]}
= \int\left(\int f[x,\cdot] g[\cdot,z] d\beta^x_{[x,y]}\right)d\lambda^{[x]}([x,y])
= \int f[x,y]g[y,z] d\alpha^{r(x)}(y).
\end{gather*}
From this expression of the convolution product, one can see that the condition (iii) of Def\/ini\-tion~\ref{Haar} is
satisf\/ied.
The adjoint property (ii) is also easily checked.
\end{proof}

Let us formalize a~def\/inition used earlier.

\begin{Definition}
A~locally compact hypergroupoid with Haar system of the form $(X*X)/G$ given by the theorem, where~$G$ is a~locally
compact groupoid with Haar system and~$X$ is a~proper~$G$-space endowed with a~continuous~$G$-equivariant system of
measures is called a~spatial hypergroupoid.
\end{Definition}

Given a~locally compact hypergroupoid with Haar system $(H,\lambda)$, we def\/ine the convolution product of~$f$ and~$g$
in $C_c(H)$ by
\begin{gather*}
f*g(x)=\int f(x*y)g(y^*)d\lambda^{s(x)}(y)=\int f(y)g(y^**x)d\lambda^{r(x)}(y)
\end{gather*}
and the involution by $f^*(x)=\overline{f(x^*)}$.

\begin{Proposition}
Let $(H,\lambda)$ be a~locally compact hypergroupoid with Haar system.
Endowed with the convolution product, the involution and the inductive limit topology, $C_c(H)$ is a~topological
$*$-algebra.
\end{Proposition}

The proof is very much like the groupoid case (see~\cite[Proposition~II.1.1]{ren:approach}) and will not be given here.

The def\/inition of the full $C^*$-algebra $C^*(H,\lambda)$ and the reduced $C^*$-algebra $C^*_r(G,\lambda)$ follows the
same line as the groupoid case.

\begin{Definition}
The full $C^*$-algebra $C^*(H,\lambda)$ of a~locally compact hypergroupoid with Haar system $(H,\lambda)$ is the
completion of the $*$-algebra $C_c(H)$ with respect to the full norm
\begin{gather*}
\|f\|=\sup\{\|L(f)\|: L~\text{non-degenerate and I-bounded representation of}~C_c(H)\},
\end{gather*}
where the I-norm of $f\in C_c(H)$ is
\begin{gather*}
\|f\|_I=\max\left(\sup_{u\in H^{(0)}}\int |f| d\lambda^u, \sup_{u\in H^{(0)}}\int |f^*| d\lambda^u\right).
\end{gather*}
\end{Definition}

In the case of a~hypergroup (i.e.~$H^{(0)}$ has only one element), this def\/inition agrees with that given
in~\cite{her:1992, her:1995}, except that we def\/ine the involution without using the modular function.

\begin{Definition}
The reduced $C^*$-algebra $C_r^*(H,\lambda)$ of a~locally compact hypergroupoid with Haar system $(H,\lambda)$ is the
completion of $C_c(H)$ for the reduced norm
\begin{gather*}
\|f\|_r=\sup\big\{\|L_u(f)\|: u\in H^{(0)}\big\},
\end{gather*}
where $L_u$ is the $*$-representation of the $*$-algebra $C_c(H)$ on the Hilbert space $L^2(H_u,\lambda_u)$ def\/ined by
$\displaystyle
L_u(f)\xi(x)=\int f(x*y) \xi(y^*)d\lambda^u(y)$ for $f\in C_c(H)$ and $\xi\in C_c(H_u)$.
\end{Definition}

When $H=(X*X)/G$ is a~spatial hypergroupoid, where~$G$ is a~locally compact groupoid with Haar system~$\lambda$ and~$X$
is a~proper~$G$-space with equivariant sytem~$\alpha$, we have def\/ined the $C^*$-algebra $C_G^*(H)$ by considering only
the representations of $C_c(H)$ induced by representations of~$G$.
Since they are I-bounded, $C_G^*(H)$ is a~quotient of $C^*(H)$.
The examples of~\cite[Subsection~15.5]{jew:hyper} or~\cite[Section~4]{her:1995}) show that $C_G^*(H)$ may be a~strict quotient of
$C^*(H)$.
On the other hand, the completion of $C_c(H)$ obtained from the regular representation of~$G$ agrees with the reduced
$C^*$-algebra $C_r^*(H)$.

\subsection*{Acknowledgements}

I thank S.~Echterhof\/f for the references~\cite{her:1992, her:1995} and J.~Brown, R.D.~Holkar, G.~Nagy and S.~Reznikof\/f
for fruitful discussions.
I~thank the referees for their great help to improve the presentation of the manuscript.

\pdfbookmark[1]{References}{ref}
\LastPageEnding


\begin{thebibliography}{99}
\footnotesize\itemsep=1.5pt

\bibitem{adr:amenable}
Anantharaman-Delaroche C., Renault J., Amenable groupoids, \textit{Monographies
  de L'Enseignement Math\'ematique}, Vol.~36, L'Enseignement Math\'ematique, Geneva, 2000.

\bibitem{bh:hyper}
Bloom W.R., Heyer H., Harmonic analysis of probability measures on hypergroups,
  \href{http://dx.doi.org/10.1515/9783110877595}{\textit{de Gruyter Studies in Mathematics}}, Vol.~20, Walter de Gruyter \&
  Co., Berlin, 1995.

\bibitem{bou:int7-8}
Bourbaki N., \'{E}l\'ements de math\'ematique. {F}ascicule~{XXIX}.
  {L}ivre~{VI}: {I}nt\'egration. {C}hapitre~7: {M}esure de {H}aar.
  {C}hapitre~8: {C}onvolution et repr\'esentations, \textit{Actualit\'es Scientifiques
  et Industrielles}, No.~1306, Hermann, Paris, 1963.

\bibitem{bs:morphisms}
Buneci M., Stachura P., Morphisms of locally compact groupoids endowed with
  Haar systems, \href{http://arxiv.org/abs/math.OA/0511613}{math.OA/0511613}.

\bibitem{dun:hyper}
Dunkl C.F., The measure algebra of a locally compact hypergroup, \href{http://dx.doi.org/10.1090/S0002-9947-1973-0320635-2}{\textit{Trans.
  Amer. Math. Soc.}} \textbf{179} (1973), 331--348.

\bibitem{hkk:hyper}
Hauenschild W., Kaniuth E., Kumar A., Harmonic analysis on central hypergroups
  and induced representations, \href{http://dx.doi.org/10.2140/pjm.1984.110.83}{\textit{Pacific~J. Math.}} \textbf{110} (1984),
  83--112.

\bibitem{her:1992}
Hermann P., Induced representations of hypergroups, \href{http://dx.doi.org/10.1007/BF02571455}{\textit{Math.~Z.}}
  \textbf{211} (1992), 687--699.

\bibitem{her:1995}
Hermann P., Representations of double coset hypergroups and induced
  representations, \href{http://dx.doi.org/10.1007/BF02567801}{\textit{Manuscripta Math.}} \textbf{88} (1995), 1--24.

\bibitem{hs:morphisms}
Hilsum M., Skandalis G., Morphismes {$K$}-orient\'es d'espaces de feuilles et
  fonctorialit\'e en th\'eorie de {K}asparov (d'apr\`es une conjecture d'{A}.
  {C}onnes), \textit{Ann. Sci. \'Ecole Norm. Sup.~(4)} \textbf{20} (1987),
  325--390.

\bibitem{hol:thesis}
Holkar R.D., Topological construction of $C^{\ast}$-correspondences for
  groupoid $C^{\ast}$-algebras, Ph.D.~thesis, G\"ottingen University, in
  preparation.

\bibitem{hr:2013}
Holkar R.D., Renault J., Hypergroupoids and {$C^*$}-algebras,
  \href{http://dx.doi.org/10.1016/j.crma.2013.11.003}{\textit{C.~R.~Math. Acad. Sci. Paris}} \textbf{351} (2013), 911--914,
  \href{http://arxiv.org/abs/1403.3424}{arXiv:1403.3424}.

\bibitem{jew:hyper}
Jewett R.I., Spaces with an abstract convolution of measures, \href{http://dx.doi.org/10.1016/0001-8708(75)90002-X}{\textit{Adv.
  Math.}} \textbf{18} (1975), 1--101.

\bibitem{kpc:hyper}
Kalyuzhnyi A.A., Podkolzin G.B., Chapovsky Yu.A., Harmonic analysis on a locally
  compact hypergroup, \textit{Methods Funct. Anal. Topology} \textbf{16}
  (2010), 304--332.

\bibitem{lan:2001}
Landsman N.P., Operator algebras and {P}oisson manifolds associated to
  groupoids, \href{http://dx.doi.org/10.1007/s002200100496}{\textit{Comm. Math. Phys.}} \textbf{222} (2001), 97--116.

\bibitem{mm:1999}
Macho~Stadler M., O'uchi M., Correspondence of groupoid {$C^\ast$}-algebras,
  \textit{J.~Operator Theory} \textbf{42} (1999), 103--119.

\bibitem{mac:induced}
Mackey G.W., Induced representations of locally compact groups.~{I},
  \textit{Ann. of Math.} \textbf{55} (1952), 101--139.

\bibitem{mrc:1999}
Mr{\v{c}}un J., Functoriality of the bimodule associated to a
  {H}ilsum--{S}kandalis map, \href{http://dx.doi.org/10.1023/A:1007773511327}{\textit{$K$-Theory}} \textbf{18} (1999), 235--253.

\bibitem{mrw:Morita}
Muhly P.S., Renault J., Williams D.P., Equivalence and isomorphism for
  groupoid {$C^\ast$}-algebras, \textit{J.~Operator Theory} \textbf{17} (1987),
  3--22.

\bibitem{ram:topologies}
Ramsay A., Topologies on measured groupoids, \href{http://dx.doi.org/10.1016/0022-1236(82)90110-0}{\textit{J.~Funct. Anal.}}
  \textbf{47} (1982), 314--343.

\bibitem{ram:dichotomy}
Ramsay A., The {M}ackey--{G}limm dichotomy for foliations and other {P}olish
  groupoids, \href{http://dx.doi.org/10.1016/0022-1236(90)90018-G}{\textit{J.~Funct. Anal.}} \textbf{94} (1990), 358--374.

\bibitem{ren:approach}
Renault J., A groupoid approach to {$C^{\ast}$}-algebras, \textit{Lecture Notes
  in Mathematics}, Vol.~793, Springer, Berlin, 1980.

\bibitem{ren:halifax}
Renault J., Strong Morita equivalence of groupoid $C^*$-algebras, Talk given in the Operator Algebras Section of the Annual Summer Meeting of the Canadian Mathematical Society (Halifax, May 27--30, 1981).


\bibitem{ren:kingston}
Renault J., {$C^{\ast} $}-algebras of groupoids and foliations, in Operator
  Algebras and Applications, {P}art {I} ({K}ingston, {O}nt., 1980),
  \textit{Proc. Sympos. Pure Math.}, Vol.~38, Amer. Math. Soc., Providence,
  R.I., 1982, 339--350.

\bibitem{ren:representations}
Renault J., Repr\'esentation des produits crois\'es d'alg\`ebres de groupo\"\i
  des, \textit{J.~Operator Theory} \textbf{18} (1987), 67--97.

\bibitem{ren:ideal}
Renault J., The ideal structure of groupoid crossed product
  {$C^\ast$}-algebras, \textit{J.~Operator Theory} \textbf{25} (1991), 3--36.

\bibitem{rie:induced}
Rief\/fel M.A., Induced representations of {$C^{\ast} $}-algebras, \href{http://dx.doi.org/10.1016/0001-8708(74)90068-1}{\textit{Adv.
  Math.}} \textbf{13} (1974), 176--257.

\bibitem{rie:kingston}
Rief\/fel M.A., Applications of strong {M}orita equivalence to transformation
  group {$C^{\ast}$}-algebras, in Operator Algebras and Applications, {P}art
  {I} ({K}ingston, {O}nt., 1980), \textit{Proc. Sympos. Pure Math.}, Vol.~38,
  Amer. Math. Soc., Providence, R.I., 1982, 299--310.

\bibitem{spe:hyper}
Spector R., Aper\c cu de la th\'eorie des hypergroupes, in Analyse harmonique
  sur les groupes de {L}ie ({S}\'em. {N}ancy-{S}trasbourg, 1973-75),
  \textit{Lecture Notes in Math.}, Vol. 497, Springer, Berlin, 1975, 643--673.

\bibitem{tu:2004}
Tu J.-L., Non-{H}ausdorf\/f groupoids, proper actions and {$K$}-theory,
  \textit{Doc. Math.} \textbf{9} (2004), 565--597, \href{http://arxiv.org/abs/math.OA/0403071}{math.OA/0403071}.

\end{thebibliography}
\end{document}